\documentclass[a4paper,11pt]{amsart}
\usepackage{amsmath,amsfonts,amssymb,amsthm}
\usepackage{graphics}
\usepackage{epsfig}
\usepackage{esint}
\newcommand{\lessim}{\stackrel{<}{\sim}}

\newcommand{\ignore}[1]{}

\newtheorem{theorem}{Theorem}

\newtheorem{lemma}{Lemma}
\newtheorem{corollary}{Corollary}

\newcommand{\R}{\mathbb{R}}

\newcommand{\Z}{\mathbb{Z}}
\newcommand{\eps}{\varepsilon}

\usepackage{color}
\definecolor{darkred}{rgb}{0.9,0.1,0.1}
\begin{document}

\title{Quasilinear SPDEs in divergence-form}
\author{Hendrik Weber and Felix Otto}

\begin{abstract}
We develop a solution theory  in H\"older spaces 
for a  quasilinear stochastic PDE driven by an additive noise.
The  key ingredients are two deterministic PDE Lemmas which 
establish a priori H\"older bounds for an equation with irregular
right hand side written in divergence form. We apply these 
deterministic bounds to the case of a noise term which is 
white in time and trace class in space to obtain
stretched exponential bounds 
for the H\"older semi-norms of the solution for the stochastic equation.

\end{abstract}

\maketitle

\section{Introduction}
We are interested in the quasi-linear equation
\begin{align}\label{i1}
\partial_t u - \nabla \cdot A(\nabla u) = \xi, 
\end{align}
with unknown $u \colon \R_t \times \R^d_x \to \R$
and for $A \colon \R^d \to \R^d$ which is uniformly elliptic, see \eqref{h01} and \eqref{h35} below for precise assumptions. 
The right hand side $\xi$ 
represents an irregular Gaussian  stochastic noise term which is white in time and coloured in 
space. We show the existence and uniqueness of 
solutions to \eqref{i1} as well as a stretched exponential moment bound on a space-time H\"older 
semi-norm for $\nabla u$ under a suitable condition on the spatial covariance operator for $\xi$.

\medskip

Quasi-linear stochastic PDE such as \eqref{i1} have been treated since the 70s, important contributions 
include \cite{pardoux1975equations,krylov1981stochastic} see also \cite{prevot2007concise}
for a more recent presentation. These works rely on the theory of monotone operators and 
yield solutions which  satisfy 
\[
\sup_{0 \leq t \leq T} \int u^2(t,x) dx + \int_0^T \int | \nabla u (t,x) |^2 dx dt < \infty
\]
for all $T<\infty$ almost surely. In fact, these methods allow for much more general equations: generalisations include
on the one hand more general non-linear operators, e.g.  degenerate cases such as the porous medium equation, 
but in particular also the case of multiplicative noise, i.e. the right hand side $\xi$ is replaced by 
$\sigma(u) \xi$ and the time integral of this product is interpreted as a stochastic integral, the latter requiring 
to introduce stochastic machinery such as filtrations, adapted processes etc.

\medskip

In this article we restrict ourselves to the case of non-degenerate $A$ and additive noise and 
develop a  regularity theory in H\"older spaces. Restricting ourself to additive noise
enables us to avoid stochastic integrals and use  purely deterministic arguments. 
The challenge then becomes to develop a theory for \eqref{i1} where the right hand side is only controlled 
in a low regularity norm.
%
%
%One key issue when deriving regularity bounds if the right hand side has very much more spatial regularity 
%than time regularity. In fact, $\int_0^t \xi(t,x) dt $ can be defined and in our case it is a well behaved $C^{\alpha}$
%function, but the time regularity is bad the time derivative of Brownian motion $W$.
%
%
%
 For this we will only access $\xi$ through the solution of the linear stochastic heat equation
\begin{align}\label{i1a}
\partial_t v -\Delta v = \xi,
\end{align}
and rewrite \eqref{i1} as 
\begin{align}\label{i1b}
\partial_t u - \nabla \cdot A(\nabla u)  = \partial_t v -\Delta v .
\end{align}
To stress the divergence form of  the right hand side, we relabel those terms 
and write
\begin{align}\label{i2}
\partial_t u - \nabla \cdot A(\nabla u) = \partial_t v  -\nabla \cdot j ,
\end{align}
for  $j = \nabla g$.
In Lemmas~\ref{L1} and \ref{L2} below we establish optimal interior a priori H\"older bounds for 
$\nabla u$ in \eqref{i2} in  terms of  the parabolic space-time H\"older norms 
$[\nabla v]_\alpha$ and $[j]_\alpha$ (parabolic H\"older norms are defined in \eqref{h40} and \eqref{h36} below), 
first in Lemma~\ref{L1} for a small $\alpha_0$ 
using the celebrated De Giorgi-Nash Theorem, and then in Lemma~\ref{L2} for arbitrary $\alpha$. 
We do not assume time regularity for $v$ and thus, although we have the optimal $\alpha$ regularity for
$\nabla u$, we do not get the matching $\frac{1+\alpha}{2}$ temporal regularity for $u$. This corresponds
exactly the fact that in our stochastic application the right hand side $\xi$ is white in time and we cannot expect to 
get more than the Brownian temporal regularity $\frac12-$  for $u$. It turns out however, that $ u-v$ 
is better behaved and we are able to control the full parabolic $1+\alpha$ H\"older semi-norm $[u-v]_{1+\alpha}$.

\medskip
Our main result, Theorem~\ref{t:main},  illustrates an application of these bounds to construct solutions 
in the random case. In order to avoid having to deal with the large scale behaviour of solutions 
we restrict ourselves to the simplest possible setting and impose that the noise $\xi$  is $1$-periodic in 
all spatial direction and additionally compactly supported in time, say on the interval $t \in [0,1]$;
we then construct solutions $u$ which satisfy $u_{|t\leq 0} =0$.

\begin{theorem}\label{t:main}
Let $A$ be uniformly elliptic with ellipticity contrast $\lambda$ and let $D A$ be Lipschitz continuous  with constant $\Lambda$ (in the sense of   \eqref{h01} and \eqref{h35} below). 
Let $\alpha_0 = \alpha_0(d,\lambda) $ be as in Lemma~\ref{L1}. 
Let $v$ be given by \eqref{e:def:g} for a covariance operator $K$ satisfying \eqref{condition} for some $s>d$.
Then for almost all realisations of $v$, there exists a unique $u = u(t,x)$ with the following properties:
\begin{itemize}
\item $u$ is  continuous, $1$-periodic in all spatial directions (i.e. $u(t,x) = u(t,x+k)$ for all $k \in \Z^d$) and $u_{|t \leq 0} = 0$.
\item $[\nabla u ]_\alpha, [u-v]_{1+\alpha} < \infty $ for $\alpha < \min \{ \frac{s-d}{2},1\}$.
\item $u$ solves \eqref{i1b} in the distributional sense, i.e. for all Schwartz functions $\varphi \in \mathcal{S}(\R \times \R^d)$
\begin{align*}
&- \int \int  \partial_t \varphi u dx dt+   \int \int  \nabla \varphi \cdot A(\nabla u) dx dt \\
& \qquad = - \int \int  \partial_t \varphi v dx dt + \int \int   \nabla \varphi  \cdot \nabla v \; dx dt  .
\end{align*}
\end{itemize}
Furthermore, for $\alpha < \min\{ \frac{s-d}{2}, 1\}$ there exists $C= C(d,\lambda,\Lambda,\alpha,s)<\infty$ such that 
\begin{align}\label{sem}
\Big\langle \exp \Big( \frac{1}{C} \big([\nabla u]_\alpha + [u-v]_{1+\alpha}\big)_\alpha^{2 \min\{ 1, \frac{\alpha_0}{\alpha} \}}   \Big) \Big\rangle < \infty,
\end{align}
where $\langle \cdot \rangle$ denotes the expectation with respect to the probability distribution of $v$.
\end{theorem}
Several higher regularity results for quasilinear SPDE were derived in the last years: both \cite{denis2005p} and
\cite{dareiotis2017local}   considered parabolic equations with a uniformly elliptic leading term with only 
measurable coefficients and a gradient-dependent noise coefficient. They derived stochastic $L^p$ bounds
on the space-time $L^\infty$ norm of solutions as well as a stochastic Harnack inequality. In 
\cite{gess2014random} a stochastic porous medium equation driven by a multiplicative noise of the form 
$\sum_{k=1}^N f_k u_k \circ d\beta^k$ was analysed and using a transformation which removes the noise term of this 
particular structure, uniform continuity of solutions was shown.
The recent work
  \cite{debussche2015regularity},  where a H\"older theory for the quasilinear
  stochastic PDE
\begin{equation}\label{dh}
\partial_t u - \nabla \cdot (A(u) \nabla u) = H(u) \xi
\end{equation}
is developed, is probably closest to the analysis presented here.
The key idea of their analysis is to consider first the auxiliary equation
\begin{align*}
\partial_t z - \Delta z=  H(u) \xi, 
\end{align*}
and to then use 
the De Giorgi-Nash Theorem to get an a priori bound on the remainder
$w = u-z$. We pursue a somewhat similar strategy, and work with the 
equation for $w= u-v$. However, \eqref{i1} is more non-linear than  \eqref{dh}
and we cannot work with the equation for $w$ directly. Instead, we linearise it by 
taking spatial differences, see \eqref{h06} below.

\medskip

In a previous version of this work, see \cite{OttoWeber}, we treated a quasilinear equation
\begin{align}\label{i10}
\frac{1}{T} u  +\partial_t u - \partial_x^2 \pi (u) = \bar{\xi},
\end{align}
where $\bar{\xi}$ is a space-time white noise over $\R_t \times \R_x$, and derived a stretched exponential 
moment bound akin to \eqref{sem} on the H\"older semi-norms $[u]_\alpha$. 
The results in the present article 
contain this result, up to the different treatment of large scales. Indeed, specialising \eqref{i1} to the case $d=1$ and differentiating 
with respect to $x$ yields for $\bar{u} = \partial_x u$ 
\begin{align*}
\partial_t \bar{u} - \partial_x^2 A(\bar{u}) = \partial_x \xi,
\end{align*}
which coincides with \eqref{i10}, noting that our assumptions on $\xi$ cover the case where 
$\bar{\xi}=\partial_x \xi$ is a space-time white noise in one spatial dimension, and that in the one-dimensional case our 
assumptions on $A$ coincide with the assumptions imposed on $\pi $ in \cite{OttoWeber}.

\section{Setting}

We are interested in the quasi-linear parabolic equation
\begin{align}\label{h02}
\partial_tu-\nabla\cdot A(\nabla u)=\partial_tv+\nabla\cdot j
\end{align}
with a rough right hand side as described by $v$ and $j$. We present a Schauder-theory where we estimate
the solution $\nabla u$ by the data $(v,j)$ in the H\"older space $C^\alpha$,
always with respect to to the parabolic distance 
\begin{align}\label{h37}
d((t',x'),(t,x)):=\sqrt{|t'-t|}+|x'-x|.
\end{align}
This is slightly different from the standard Schauder theory in $C^{1,\alpha}$, which cannot
be applied due to the right hand side term $\partial_tv$ that is irregular in time. In fact we shall
control the $C^{1,\alpha}$-semi norm of $w:=u-v$
\begin{align}\label{h40}
[w]_{1+\alpha}:=[\nabla w]_{\alpha}+\sup_{t,t',x}\frac{|w(t,x)-w(t',x)|}{\sqrt{|t-t'|}^{\alpha+1}},
\end{align}
with $[\cdot]_\alpha$ defined in (\ref{h36}).
We make two assumptions on the nonlinearity $A\colon\mathbb{R}^d\rightarrow\mathbb{R}^d$ 
in form of assumptions on the tensor field given by the derivative matrix $DA$:
\begin{itemize}
\item $DA$ is uniformly elliptic in the sense that there exists a constant $\lambda>0$ such that
\begin{align}\label{h01}
\xi\cdot DA(q)\xi\ge\lambda|\xi|^2\quad\mbox{and}\quad|DA(q)\xi|\le|\xi|
\quad\mbox{for all vectors}\;q,\xi.
\end{align}
Here, without loss of generality  we normalized the upper bound to unity.
We will make use of (\ref{h01}) in the following form: For every spatial shift vector $y\in\mathbb{R}^d$
we will work with the increment operator $\delta_yu(t,x)=u(t,x+y)-u(t,x)$ and use the
chain-type rule
\begin{align}\label{h03}
\delta_yA(\nabla u)=a_y\delta_y\nabla u
\end{align}
where
\begin{align*}
a_y(t,x)=\int_0^1DA(\theta\nabla u(t,x+y)+(1-\theta)\nabla u(t,x))d\theta.
\end{align*}
Then (\ref{h01}) ensures that for all $y$ we have uniform ellipticity of $a_y$:
\begin{align}\label{h07}
\xi\cdot a_y(t,x)\xi\ge\lambda|\xi|^2\quad\mbox{and}\quad|a_y(t,x)\xi|\le|\xi|
\quad\mbox{for all}\;(t,x),\xi.
\end{align}
\item $DA$ is globally Lipschitz in the sense that there exists a constant $\Lambda<\infty$ such that
\begin{align}\label{h35}
|DA(q')-DA(q)|\le \Lambda|q'-q|\quad\mbox{for all}\;q,q'.
\end{align}
We will make use of (\ref{h35}) in the following form: For any exponent $\beta\in(0,1]$ we have
the following estimate on the level of H\"older norms 
\begin{align}\label{h20}
[a_y]_{\beta}\le \Lambda [\nabla u]_{\beta}\quad\mbox{and}\quad[A(\nabla u)]_\beta\le\Lambda[\nabla u]_\beta,
\end{align}
where $[\cdot]_\beta$ denotes the (parabolic) H\"older semi-norm on space-time $\mathbb{R}_t\times\mathbb{R}_x^d$
\begin{align}\label{h36}
[a]_\beta:=\sup_{z'\not=z\in\mathbb{R}\times\mathbb{R}^d}\frac{|a(z')-a(z)|}{d^\beta(z',z)}
\end{align}
and $d$ the parabolic distance, cf (\ref{h37}).
%
%\begin{align}\label{h37}
%d((t',x'),(t,x)):=\sqrt{|t'-t|}+|x'-x|.
%\end{align}
\end{itemize}

The use equation (\ref{h02}) exclusively in the following form: We apply the increment operator 
$\delta_y$ to it and obtain by (\ref{h03})
\begin{align*}
\partial_t\delta_yu-\nabla\cdot a_y\nabla\delta_y u=\partial_t\delta_yv+\nabla\cdot \delta_yj,
\end{align*}
which in terms of the difference $w:=u-v$ we rewrite as
\begin{align}\label{h06}
\partial_t\delta_yw-\nabla\cdot a_y\nabla\delta_y w=\nabla\cdot(a_y\nabla\delta_yv+\delta_yj).
\end{align}

We establish our form of $C^{1,\alpha}$-Schauder theory, cf Corollary \ref{main}, in two lemmas. While 
the Lemma \ref{L1} just relies on the uniform ellipticity (\ref{h01}) and crucially
uses the $C^\alpha$-a priori estimate of De Giorgi and Nash,  
Lemma \ref{L2} uses also the Lipschitz continuity (\ref{h35}) and proceeds by a Schauder-type
argument.

\begin{lemma}\label{L1}
There exists an exponent $\alpha_1=\alpha_1(d,\lambda)\in(0,1)$ such that
for any exponent $\alpha_0\in(0,\alpha_1)$ we have
\begin{align}\label{h04}
[\nabla u]_{\alpha_0}\le C(d,\lambda,\alpha_0)\big([\nabla v]_{\alpha_0}+[j]_{\alpha_0}\big)
\end{align}
provided we already have the qualitative information that the left hand side is finite.
\end{lemma}

\begin{lemma}\label{L2}
Let $\alpha_0$ be as in Lemma \ref{L1} and suppose that $L$ is so small that
\begin{align}\label{h19}
[\nabla u]_{\alpha_0,P_{2L}}\le L^{-\alpha_0},
\end{align}
where $P_R:=(-R^2,0)\times B_R$ denotes the (centered) parabolic cylinder of size $R$
and $[\cdot]_{\beta,P_R}$ the $\beta$-H\"older semi-norm restricted to this set.
Then we have for any exponent $\alpha\in[\alpha_0,1)$
\begin{align}\label{h05}
[\nabla u]_{\alpha,P_L}\le 
C(d,\lambda,\Lambda,\alpha_0,\alpha)\big(L^{-\alpha}+[\nabla v]_{\alpha,P_{2L}}+[j]_{\alpha,P_{2L}}\big).
\end{align}
\end{lemma}

\begin{corollary}\label{main}
Let $\alpha_0$ be as in Lemma \ref{L1}. Then we have for any exponent $\alpha\in(0,1)$
\begin{align*}%\label{h05}
\lefteqn{[u-v]_{1+\alpha}+[\nabla u]_{\alpha}}\nonumber\\
&\le C(d,\lambda,\Lambda,\alpha)
\Big(\big([\nabla v]_{\alpha_0}+[j]_{\alpha_0}\big)^\frac{\alpha}{\alpha_0}+
     \big([\nabla v]_{\alpha}+[j]_{\alpha}\big)\Big),
\end{align*}
provided we already have the qualitative information that $[\nabla u]_{\alpha_0}<\infty$.
\end{corollary}

%%%%%%%%%%%%%%%%%%%%%%%%%%%%%%%%%%%%%%%%%%%%%%%%%%%%%%%%%%%%%%%%%%%%%%%%%%%%%%%%%%%%%

The example that we have a mind is the case where the right hand side is 
a stochastic noise which is white in time but coloured in space. Such a noise term 
is described by a Gaussian random distribution $\xi$ over $(t,x) \in \R \times \R^d$, whose probability distribution is characterized 
by having zero mean and
\begin{align*}
\Big\langle  \Big(  \int \xi(t,x) \varphi(t,x)  dx dt \Big)^2 \Big\rangle = \int_0^T \int \int    K(x,y)\varphi(t,x) \varphi(t,y) dx dy dt  ,
\end{align*}
where the spatial correlation $K$ is given by the kernel of a regularising operator.
We denote by $v$  the solution of the constant-coefficient heat equation 
\begin{align}\label{eq:g}
(\partial_t - \Delta ) v = \xi.
\end{align}
Under suitable conditions on the kernel $K$ it is known that $\nabla v$ is regular enough, i.e. $\alpha$-H\"older continuous, 
to apply the deterministic theory. We illustrate this in the simplest possible case, where $\xi$ is assumed to be $1$-periodic in 
all spatial directions and in addition localised to a compact time interval, say the interval $[0,1]$. If we assume in addition that the 
the probability distribution of $\xi$ is translation invariant in the spatial directions we have the following convenient Fourier series 
representation 
\begin{align}\label{e:def:xi}
\xi(t,x) = \sum_{k \in (2 \pi \Z)^d } e^{ik \cdot x}  \sqrt{\hat{K}(k)} \mathbf{1}_{[0,1]}(t) \dot{\beta}_k(t) , 
\end{align}
where the $\beta_k$ are complex valued standard Brownian motions (i.e. real and imaginary parts are independent and  satisfy $\langle \mathfrak{R} (\beta_k(t))^2 \rangle =
\langle \mathfrak{I} (\beta_k(t))^2 \rangle = \frac{t}{\sqrt{2}}$), which are independent up to the constraint $\beta_k = \overline{\beta}_{-k}$, which assures that $\xi$ is real-valued,
and $\dot{\beta}_k(t)$ stands for the distributional time derivative. The almost sure convergence of \eqref{e:def:xi} in the space of distributions can be shown relatively 
easily, but we adopt the slightly simpler framework to only work with $v$ which we  define by its Fourier series representation
\begin{align}\label{e:def:g}
v(t,x) &= \sum_{k \in (2 \pi \Z)^d}  \sqrt{\hat{K}(k)} e^{i  k \cdot x }\int_0^{\min\{ t, 1\}} e^{-(t-s) |k|^2} d \beta_k(s).
\end{align}

In order to ensure that the gradient is well behaved we impose that there exists $s >d $ such that for  $k\in (2 \pi \Z)^d$
\begin{align}\label{condition}
\hat{K}(k) \leq  (1 + |k|^2 )^{-\frac{s}{2}},
\end{align}
where we have set the normalisation equal to $1$  without loss of generality.
 Incidentally, this condition on $s$  says precisely that the spatial covariance operator $K$
is of trace class. Then we have the following Lemma.
\begin{lemma}\label{lem:stoch}
Let $v(t,x)$ be given by \eqref{e:def:g} for $t> 0$ and $v_{t\leq 0}=0$.
Then there exists a $C = C(\alpha,s)< \infty$ such that for $\alpha < \min\{ \frac{s - d}{2}, 1 \}$ 
\begin{align*}
\Big\langle  \exp\Big (\frac{1}{C} [ \nabla v]_\alpha^2  \Big)\Big\rangle < \infty,
\end{align*}
where $\langle \cdot \rangle$ represents the expectation with respect to the 
probability distribution of $v$.
\end{lemma}

%%%%%%%%%%%%%%%%%%%%%%%%%%%%%%%%%%%%%%%%%%%%%%%%%%%%%%%%%%%%%%%%%%%%%

\section{Proof of Theorem~\ref{t:main}}
We prove a slightly stronger statement than announced in the Theorem: We assume we are given continuous functions $v$ and $j$ with
$[\nabla v]_\alpha$, $[j]_\alpha < \infty$ for an $\alpha \in (0,1)$, which are $1$-periodic in each spatial direction and with $v_{|t  \leq 0} = j_{|t \leq 0} =0$. 
We show that there exists a unique function $u$ which is one-periodic in each spatial direction, satisfies $u_{t \leq 0} =0$ and 
which satisfies 
\begin{align}
\notag
&- \int \int  \partial_t \varphi u dx dt+   \int \int  \nabla \varphi \cdot A(\nabla u) dx dt \\
\label{pt1}
& \qquad = - \int \int  \partial_t \varphi v dx dt + \int \int   \nabla \varphi  \cdot \nabla j \; dx dt  ,
\end{align}
for each Schwartz function $\varphi$. In addition, we have the bound 
\begin{align}\label{h05A}
[\nabla u]_{\alpha} + [u-v]_{1+\alpha} \lesssim N ,
\end{align}
where $\lesssim$ means $\leq C(d,\lambda,\Lambda,\alpha)$ and
 $N = \big([\nabla v]_{\alpha}+[j]_{\alpha}\big)^\frac{\alpha}{\alpha_0}+
     \big([\nabla v]_{\alpha}+[j]_{\alpha}\big)$.
The desired statement then follows, by applying this to the case where $v$ is given by \eqref{e:def:g}, $j = \nabla v$ and invoking Lemma~\ref{lem:stoch}.
From now on, all functions $u,v,w$ etc. appearing in the proof are assumed to be one-periodic in all space directions. 
Under this periodicity assumption the weak formulation \eqref{pt1} can be restated equivalently by 
replacing the space integrals over $\R^d$ by integrals over $[0,1]^d$ and assuming that the test functions $\varphi$ 
are also periodic. 

\medskip

The existence of solutions follows by approximation through regularisation. Let $j_\eps$,
$v_\eps  $ be space-time regularisations of $j$, $v$ satisfying $[j_\eps]_\alpha \leq [j]_\alpha$, $[\nabla v_\eps]_\alpha \leq [\nabla v]_\alpha  $
and such that ${v_{\eps}}_{|t\leq -\eps} = {j_{\eps}}_{|t\leq -\eps} = 0$.
%To be more precise, we let $\eta \colon \R \times \R^d$ be smooth, non-negative, compactly supported in $\{(t,x) \colon |t| + |x| \leq 1 \}$ with
%$\int\int \eta dt dx = 1$, then let $\eta_\eps (t,x) = \eps^{-d-1} \eta(\frac{t}{\eps}, \frac{x}{\eps})$  and define $j_\eps = j \ast \eta_\eps$,
%$v_\eps = v \ast \eta_\eps$ as the space-time convolutions of $j$ and $v$ with $\eta_\eps$ (in fact here we could also 
%choose to define $\eta_\eps$ with a parabolic rescaling, but our choice will be more convenient below).
 Then by classical theory  
 there exists a unique classical solution $u_\eps$ for
\begin{align*}
\partial_t u_\eps - \nabla \cdot A(\nabla u_\eps) = \partial_t v_\eps - \nabla \cdot j_\eps,  \qquad {u_\eps}_{|t \leq -\eps }=0,
\end{align*}
which is one-periodic in all spatial directions (see e.g. \cite[Thm. 12.14]{Lieberman}
 for a proof in the case of Dirichlet data on a bounded spatial domain. The case of the torus is only simpler).
In this situation Corollary~\ref{main} applies and yields
\begin{align}
%\notag
[\nabla u_\eps]_{\alpha} + [u_\eps-v_\eps]_{1+\alpha} &\lesssim 
 \Big(\big([\nabla v_\eps]_{\alpha}+[j_\eps]_{\alpha}\big)^\frac{\alpha}{\alpha_0}+
     \big([\nabla v_\eps]_{\alpha}+[j_\eps]_{\alpha}\big)\Big)
     \label{m10}
\lesssim 
N.
\end{align}
This estimate together with the initial datum ${u_\eps}_{|t=-\eps} = {v_\eps}_{|t=-\eps} =0$
yields enough compactness to conclude that up to choosing a subsequence
$u_\eps -v_\eps \to w$, $\nabla (u_\eps - v_\eps) \to  \nabla w$, $\nabla u_\eps \to \nabla u$  locally uniformly
for functions $u,w$ with $u=w+v$. Furthermore, $w$ solves
\begin{align*}
\partial_t w - \nabla \cdot A(\nabla u) = \nabla j
\end{align*}
in the distributional sense. Setting $u = w +v$ we obtain \eqref{pt1} and 
the estimate \eqref{h05A} follows by passing to the limit in \eqref{m10} using lower semi-continuity.

\medskip

It only remains to argue for uniqueness. Assume thus that $u^1$ and $u^2$ are one-periodic in space,
 satisfy \eqref{pt1}  and vanish for $t \leq 0$. Thus the difference $\delta u = u^1 - u^2$   satisfies 
\begin{align}\label{m25}
\partial_t \delta u = \nabla \cdot \big(A(\nabla u_1) - A(\nabla u_2) \big),
\end{align}
in the distributional sense
and $\delta u_{|t=0} =0$.
In order to show that $\delta u=0$ we aim to test \eqref{m25} equation against $\delta u$ to obtain the identity 
\begin{align}\label{pt4}
\frac{1}{2}  \int_{[0,1]^d} \delta u^2(T,x) dx = - \int_0^T \int_{[0,1]^d} \nabla \delta u \cdot \big(A(\nabla u_1) - A(\nabla u_2) \big) dx dt  ,
\end{align}
for all $T \geq 0$. Once the identity \eqref{pt4} is justified,  we can invoke the uniform ellipticity \eqref{h07} once more and obtain
the point-wise identity
\begin{align*}
\nabla \delta u \cdot \big(A(\nabla u_1) - A(\nabla u_2) \big) &= \nabla \delta u   \Big(\int_0^1 D A(\lambda \nabla u_1 +(1-\lambda) \nabla u_2) d \lambda \Big) \nabla \delta u \\
&\geq 0
\end{align*}
so that \eqref{pt4} yields $\delta u = 0$.

\medskip
It thus remains to justify \eqref{pt4}. For this we convolve \eqref{pt4} with a temporal regularising kernel at scale $\eps$ and 
then test against $\delta u_\eps$, the temporally regularised version of $\delta u$. This yields for any $T>0$
\begin{align}\label{pt4r}
\frac{1}{2}  \int_{[0,1]^d} \delta u_\eps^2(T,x) dx = - \int_0^T \int_{[0,1]^d} \nabla \delta u_\eps \cdot \big(A(\nabla u_1) - A(\nabla u_2) \big)_\eps dx dt  .
\end{align}
We can pass to the limit $\eps \to 0$ on both sides using the fact that $\delta u= (u^1- v)- (u^2-v)$ is $\frac{1+\alpha}{2}$ -H\"older in time and 
using the fact that $\nabla u^1$ and $\nabla u^2$ are $\frac{\alpha}{2}$ H\"older in time.

%%%%%%%%%%%%%%%%%%%%%%%%%%%%%%%%%%%%%%%%%%%%%%%%%%%%%%%%%%%%%%%%%%%%%%%%%%%%%%%%%%%%%

\section{Proof of Lemma \ref{L1}}

Based on (\ref{h06}) and (\ref{h07}) we have by a localized version of the H\"older a priori
estimate of De Giorgi and Nash that there exists
an exponent $\alpha_1=\alpha_1(d,\lambda)\in(0,1)$ such that for all shift vectors $y$, all
length scales $\ell$ and all space-time points $z$
\begin{align}\label{h08}
\lefteqn{[\delta_yw]_{\alpha_1,P_\ell(z)}}\nonumber\\
&\lesssim\ell^{-\alpha_1}\inf_k\|\delta_yw-k\|_{P_{2\ell}(z)}
+\ell^{1-\alpha_1}\|a_y\nabla\delta_yv+\delta_yj\|_{P_{2\ell}(z)},
\end{align}
where $\lesssim$ means $\le C(d,\lambda,\alpha_0)$, where $P_\ell(z)=(t-\ell^2,t)\times B_\ell(x)$ denotes the parabolic cylinder 
centered around $z=(t,x)$, and where $\|\cdot\|_{P_\ell(z)}$ stands for the supremum norm restricted to the set $P_\ell(z)$. 
The exponents of the $\ell$-factors in (\ref{h08}) are determined by scaling; smuggling
in the constant $k$ is possible since (\ref{h07}) is oblivious to changing $\delta_yw$ by an
additive constant. We refer to \cite[Theorem 6.28]{Lieberman} as one possible reference
(with $b\equiv0$, $c^0\equiv 0$, and $g\equiv0$ so that $k_1=\sup_{Q(R)}|f|$ in the notation
of that reference).
We fix an exponent $\alpha_0\in(0,\alpha_1)$ and take the supremum of (\ref{h08}) 
over all shift vectors $y$ with $|y|\le r$ for some $r\le\ell$ 
\begin{align}\label{h09}
\lefteqn{\sup_{|y|\le r}[\delta_yw]_{\alpha_1,P_\ell(z)}}\nonumber\\
&\lesssim\ell^{-\alpha_1}\sup_{|y|\le r}\inf_k\|\delta_yw-k\|_{P_{2\ell}(z)}
+\ell^{1-\alpha_1}\sup_{|y|\le r}\|a_y\nabla\delta_yv+\delta_yj\|_{P_{2\ell}(z)}.
\end{align}
We first estimate the right hand side terms of (\ref{h09}). We start with the second right hand side term: From the definition (\ref{h36})
of the H\"older semi-norm and that of the parabolic cylinder, we obtain
\begin{align*}
\lefteqn{\|a_y\nabla\delta_yv+\delta_yj\|_{P_{2\ell}(z)}}\nonumber\\
&\stackrel{(\ref{h07})}{\le}\|\delta_y\nabla v\|_{P_{2\ell}(z)}+\|\delta_y j\|_{P_{2\ell}(z)}
\le |y|^{\alpha_0}([\nabla v]_{\alpha_0}+[j]_{\alpha_0}),
\end{align*}
so that
\begin{align}\label{h11}
\ell^{1-\alpha_1}\sup_{|y|\le r}\|a_y\nabla\delta_yv+\delta_yj\|_{P_{2\ell}(z)}
\lesssim \ell^{1-\alpha_1}r^{\alpha_0}([\nabla v]_{\alpha_0}+[j]_{\alpha_0}).
\end{align}
We now turn to the first right hand side term of (\ref{h09}): We first note that
\begin{align}\label{h10}
\inf_k\|\delta_yw-k\|_{P_{2\ell}(z)}\le r\inf_{c}\|\nabla w-c\|_{P_{3\ell}(z)},
\end{align}
where the right hand side infimum ranges over all $c\in\mathbb{R}^d$.
Indeed, passing to $\tilde w(t,x)=w(t,x)-c\cdot y$, so that $\nabla w-c$ $=\nabla\tilde w$,
and transforming $\tilde k=k-c\cdot y$, so that $\delta_yw-k=\delta_y\tilde w-\tilde k$, 
we see that (\ref{h10}) reduces to $\|\delta_y\tilde w\|_{P_{2\ell}(z)}$ 
$\le |y|\|\nabla \tilde w\|_{P_{3\ell}(z)}$, which because of $|y|\le r\le\ell$
is a consequence of the mean-value theorem. Since obviously $\inf_{c}\|\nabla w-c\|_{P_{3\ell}(z)}$
$\le (3\ell)^{\alpha_0}[\nabla w]_{\alpha_0}$, we obtain
\begin{align}\label{h12}
\ell^{-\alpha_1}\sup_{|y|\le r}\inf_k\|\delta_yw-k\|_{P_{2\ell}(z)}
\lesssim \ell^{\alpha_0-\alpha_1}r[\nabla w]_{\alpha_0}.
\end{align}
We finally turn to the left hand side term in (\ref{h09}) and note 
\begin{align}\label{h13}
\sup_{|y|\le r}[\delta_yw]_{\alpha_1,P_\ell(z)}\ge
\sup_{|y|\le r}[\delta_yw]_{\alpha_1,P_r(z)}\ge
\frac{1}{r^{\alpha_1}}
\sup_{|y|\le r}\inf_{k}\|\delta_yw-k\|_{P_r(z)}.
\end{align}
Inserting (\ref{h11}), (\ref{h12}), and (\ref{h13}), into (\ref{h09}) we obtain
\begin{align*}%\label{h09}
\lefteqn{\frac{1}{r^{\alpha_1}}\sup_{|y|\le r}\inf_{k}\|\delta_yw-k\|_{P_r(z)}}\nonumber\\
&\lesssim \ell^{\alpha_0-\alpha_1}r[\nabla w]_{\alpha_0}
+\ell^{1-\alpha_1}r^{\alpha_0}([\nabla v]_{\alpha_0}+[j]_{\alpha_0}),
\end{align*}
which we arrange to
%Using the triangle inequality on $[\nabla w]_{\alpha_0}$ and $r\le\ell$ on the first right hand side term,
%we rearrange this to
%
\begin{align}\label{h15}
\lefteqn{\frac{1}{r^{1+\alpha_0}}\sup_{|y|\le r}\inf_{k}\|\delta_yw-k\|_{P_r(z)}}\nonumber\\
&\lesssim (\frac{r}{\ell})^{\alpha_1-\alpha_0}[\nabla w]_{\alpha_0}
+(\frac{\ell}{r})^{1-\alpha_1}([\nabla v]_{\alpha_0}+[j]_{\alpha_0}).
\end{align}

\medskip

We now argue that we are done once we establish the norm equivalence
\begin{align}\label{h14}
[\nabla w]_{\alpha_0}\lesssim\sup_{z,r}\frac{1}{r^{1+\alpha_0}}\sup_{|y|\le r}\inf_{k}\|\delta_yw-k\|_{P_r(z)}.
\end{align}
Indeed, choosing $\ell= Mr$ with $M\ge 1$ to be fixed presently, we take the supremum of (\ref{h15}) over
all radii $r$ and all space-time points $z$ to arrive at
\begin{align*}
\lefteqn{\sup_{z,r}\frac{1}{r^{1+\alpha_0}}\sup_{|y|\le r}\inf_{k}\|\delta_yw-k\|_{P_r(z)}}\nonumber\\
&\lesssim M^{\alpha_0-\alpha_1}[\nabla w]_{\alpha_0}
+M^{1-\alpha_1}([\nabla v]_{\alpha_0}+[j]_{\alpha_0}),
\end{align*}
into which we insert (\ref{h14})
\begin{align*}
[\nabla w]_{\alpha_0}\lesssim M^{\alpha_0-\alpha_1}[\nabla w]_{\alpha_0}
+M^{1-\alpha_1}([\nabla v]_{\alpha_0}+[j]_{\alpha_0}).
\end{align*}
By the triangle inequality in $[\cdot]_{\alpha_0}$ we post-process this to
\begin{align*}
[\nabla u]_{\alpha_0}\lesssim M^{\alpha_0-\alpha_1}[\nabla u]_{\alpha_0}
+M^{1-\alpha_1}([\nabla v]_{\alpha_0}+[j]_{\alpha_0}).
\end{align*}
Since by our qualitative assumption of $[\nabla u]_{\alpha_0}<\infty$, and since $\alpha_0<\alpha_1$,
we may choose $M=M(d,\lambda,\alpha_0)$ so large that this turns into the desired (\ref{h04}).

\medskip

We now turn to the norm equivalence (\ref{h14}); the elements of the argument are standard
in modern Schauder theory, in the spirit of \cite[Theorem 3.3.1]{Krylov}. 
By rotational symmetry, it is enough to establish
\begin{align}\label{h18}
[\partial_1 w]_{\alpha_0}\lesssim\sup_{z,r}\frac{1}{r^{1+\alpha_0}}\sup_{|y|\le r}\inf_{k}\|\delta_{y}w-k\|_{P_r(z)}
=:N.
\end{align}
Let $k=k(y,r,z)$ denote the optimal constant in the right hand side of (\ref{h18}).
We first argue that for arbitrary but fixed point $z$, we have for all radii $r$
\begin{align}\label{h16}
|k(2re_1,2r,z)-2k(re_1,r,z)|\lesssim Nr^{1+\alpha_0}.
\end{align}
Indeed, based on the telescoping identity $\delta_{2re_1}w$ $=\delta_{re_1}w$ $+\delta_{re_1}w(\cdot+re_1)$
we obtain by the triangle inequality the following additivity of $k$ in the $y$-variable
\begin{align*}
\lefteqn{|k(2re_1,2r,z)-2k(re_1,2r,z)|\le\|\delta_{2re_1}w-k(2re_1,2r,z)\|_{P_r(z)}}\nonumber\\
&+  \|\delta_{ re_1}w-k( re_1,2r,z)\|_{P_r(z)}
+  \|\delta_{ re_1}w(\cdot+re_1)-k( re_1,2r,z)\|_{P_r(z)}\\
&\le\|\delta_{2re_1}w-k(2re_1,2r,z)\|_{P_{2r}(z)}
+ 2\|\delta_{ re_1}w-k( re_1,2r,z)\|_{P_{2r}(z)}\\
&\le 3(2r)^{1+\alpha_0}N.
\end{align*}
Likewise, we have that $k$ only mildly depends on the $r$-variable
\begin{align*}
\lefteqn{|k(re_1,2r,z)-k(re_1,r,z)|}\nonumber\\
&\le\|\delta_{re_1}w-k(re_1,2r,z)\|_{P_{2r}(z)}+\|\delta_{re_1}w-k(re_1,r,z)\|_{P_{r}(z)}\nonumber\\
&\le 2(2r)^{1+\alpha_0}N.
\end{align*}
From the two last estimates, we obtain (\ref{h16}). Since $\alpha_0>0$, we learn from (\ref{h16}) 
that there exists a constant $c_1(z)$ such that
\begin{align*}
|\frac{1}{r}k(re_1,r,z)-c_1(z)|\lesssim Nr^{\alpha_0},
\end{align*}
along a given dyadic sequence of radii $r$.
We insert this into the definition of $N$ to obtain
\begin{align*}
\|\frac{1}{r}\delta_{re_1}w-c_1(z)\|_{P_r(z)}\lesssim Nr^{\alpha_0},
\end{align*}
from which, since in particular $u$ and thus $w$ is differentiable in the spatial variable, 
we learn that $c_1(z)=\partial_1w(z)$ so that
\begin{align}\label{h17}
\|\frac{1}{r}\delta_{re_1}w-\partial_1w(z)\|_{P_r(z)}\lesssim Nr^{\alpha_0}.
\end{align}
Since we identified the limit, this now holds for any radius $r$.
Given two points $z$, $z'$ we set $r:=2d(z,z')$, cf (\ref{h37}), and obtain
% where $d(z,z')$ $:=\sqrt{|t-t'|}+|x-x'|$ denotes the parabolic distance we have
%
\begin{align*}
\lefteqn{|\partial_1w(z)-\partial_1w(z')|}\nonumber\\
&\le\|\frac{1}{r}\delta_{re_1}w-\partial_1w(z)\|_{P_{\frac{r}{2}}(z)}
+\|\frac{1}{r}\delta_{re_1}w-\partial_1w(z')\|_{P_{\frac{r}{2}}(z)}\nonumber\\
&\le\|\frac{1}{r}\delta_{re_1}w-\partial_1w(z)\|_{P_{r}(z)}
+\|\frac{1}{r}\delta_{re_1}w-\partial_1w(z')\|_{P_{r}(z')}.
\end{align*}
Hence (\ref{h18}) follows from (\ref{h17}).

%%%%%%%%%%%%%%%%%%%%%%%%%%%%%%%%%%%%%%%%%%%%%%%%%%%%%%%%%%%%%%%%%%%%%%%%%%%%%%%%%%%%%%%%%%%%%%%%%%%%

\section{Proof of Lemma \ref{L2}}

Let the two scales $r\le \ell\le\frac{L}{4}$ be arbitrary and for the time being fixed.
Let $y$ be an arbitrary shift vector with $|y|\le r$.
By (\ref{h20}) in the localized form of $[a_y]_{\alpha_0,P_{3\ell}}$ 
$\le\Lambda[\nabla u]_{\alpha_0,P_{3\ell+r}}$ and (\ref{h19}) we have
\begin{align}\label{h21}
[a_y]_{\alpha_0,P_{3\ell}}\lesssim \ell^{-\alpha_0},
\end{align}
where $\lesssim$ stands for $\le C(d,\lambda,\Lambda,\alpha_0,\alpha)$. In conjunction with (\ref{h07})
we see that we may apply standard $C^{1,\alpha_0}$-Schauder theory to the parabolic operator 
$\partial_t-\nabla\cdot a_y\nabla$ 
when localized to $P_{3\ell}$. We see from rescaling according to $(t,x)=(\ell^2\hat t,\ell\hat x)$
that (\ref{h21}) is exactly the control on the coefficient needed so
that the constant in this localized Schauder theory is of the desired form $C(d,\lambda,\Lambda,\alpha_0,\alpha)$.
We refer to \cite[Theorem 4.8]{Lieberman} for a possible reference (with $b\equiv0$, $c\equiv0$, $g\equiv0$
in the notation of that reference).
We apply this to the increment $\delta_yw$, cf (\ref{h06}), to the effect of
\begin{align}\label{h22}
\lefteqn{\ell^{\alpha_0}[\nabla\delta_yw]_{\alpha_0,P_{2\ell}}+\|\nabla\delta_yw\|_{P_{2\ell}}}\nonumber\\
&\lesssim\ell^{-1}\inf_k\|\delta_yw-k\|_{P_{3\ell}}+\ell^{\alpha_0}[a\nabla\delta_yv+\delta_yj]_{\alpha_0,P_{3\ell}}.
\end{align}
We first argue that we may upgrade (\ref{h22}) to
\begin{align}\label{h26}
\lefteqn{\inf_c\|\nabla w-c\|_{P_r}}\nonumber\\
&\lesssim\sup_{|y|\le r}\big(\ell^{-1}\inf_k\|\delta_yw-k\|_{P_{3\ell}}
+\ell^{\alpha_0}[a\nabla\delta_yv+\delta_yj]_{\alpha_0,P_{3\ell}}\big).
\end{align}
The first ingredient in passing from (\ref{h22}) to (\ref{h26}) is the following elementary interpolation
estimate
\begin{align}\label{h27}
\inf_c\|\nabla w-c\|_{P_r}
\lesssim\sup_{|y|\le r}\big(r\|\partial_t(\delta_yw)_r\|_{P_r}+\|\nabla\delta_yw\|_{P_r}\big),
\end{align}
where $(\cdot)_r$ denotes convolution on scale $r$ in the spatial variable. Here comes the argument for (\ref{h27})
where without loss of generality  we may assume $r=1$ and restrict to estimating the first component $\partial_1w$ of the gradient. 
Given $(t,x)\in P_1$ this follows from combining
\begin{align*}
|\partial_1w(t,x)-\partial_1w(t,0)|&\le\|\nabla\delta_xw\|_{P_1},\nonumber\\
|\partial_1w(t,0)-(\delta_{e_1}w)_1(t,0)|&\lesssim\sup_{|s|\le 1}\|\nabla\delta_{se_1}w\|_{P_1},\nonumber\\
|(\delta_{e_1}w)_1(t,0)-(\delta_{e_1}w)_1(0,0)|&\le\|\partial_t(\delta_{e_1}w)_1\|_{P_1},\nonumber
\end{align*}
so that $c$ in (\ref{h27}) is given by $(\delta_{e_1}w)_1(0,0)$. The second ingredient in passing from (\ref{h22}) to
(\ref{h26}) is
\begin{align}\label{h29}
\lefteqn{r\|\partial_t(\delta_yw)_r\|_{P_r}}\nonumber\\
&\lesssim\ell^{\alpha_0}\big([\nabla\delta_yw]_{\alpha_0,P_{2\ell}}
+[a\nabla\delta_yv+\delta_yj]_{\alpha_0,P_{2\ell}}\big)+\|\nabla\delta_yw\|_{P_{2\ell}}.
\end{align}
In order to see this we apply the spatial convolution operator $(\cdot)_r$ to (\ref{h06}) to the effect of
\begin{align*}
\partial_t(\delta_yw)_r=\nabla\cdot(a_y\nabla\delta_yw+a_y\nabla\delta_yv+\delta_yj)_r.
\end{align*}
From this representation and $r\le\ell$ we obtain the estimate
\begin{align*}
\lefteqn{\|\partial_t(\delta_yw)_r\|_{P_\ell}}\nonumber\\
&\lesssim r^{\alpha_0-1}[a\nabla\delta_yw+a\nabla\delta_yv+\delta_yj]_{\alpha_0,P_{2\ell}}\nonumber\\
&\le r^{\alpha_0-1}\big([a]_{\alpha_0,P_{2\ell}}\|\nabla\delta_yw\|_{P_{2\ell}}+
\|a\|_{P_{2\ell}}[\nabla\delta_yw]_{\alpha_0,P_{2\ell}}\nonumber\\
&\quad+[a\nabla\delta_yv+\delta_yj]_{\alpha_0,P_{2\ell}}\big)\nonumber\\
&\stackrel{(\ref{h21}),(\ref{h07})}{\lesssim} 
r^{-1}(\frac{r}{\ell})^{\alpha_0}\|\nabla\delta_yw\|_{P_{2\ell}}
+r^{\alpha_0-1}\big([\nabla\delta_yw]_{\alpha_0,P_{2\ell}}
+[a\nabla\delta_yv+\delta_yj]_{\alpha_0,P_{2\ell}}\big),\nonumber
\end{align*}
which because of $r\le\ell$ yields (\ref{h29}). Inserting (\ref{h22}) into (\ref{h29}), and
the outcome into (\ref{h27}), we obtain (\ref{h26}).

\medskip

We now address the right hand side terms of (\ref{h26}). In view of (\ref{h10}) (slightly modified) we have for
the first right hand side term
\begin{align}\label{h30}
\inf_k\|\delta_yw-k\|_{P_{3\ell}}\le r\sup_{c}\|\nabla w-c\|_{P_{4\ell}}.%\le r(4\ell)^{\alpha}[\nabla w]_{P_{L}}
\end{align}
We now turn to the second right hand side term of (\ref{h26}) and note that by (\ref{h21}) and (\ref{h07})
\begin{align}\label{h23}
[a\nabla\delta_yv+\delta_yj]_{\alpha_0,P_{3\ell}}
&\le[a]_{\alpha_0,P_{3\ell}}\|\nabla\delta_yv\|_{P_{3\ell}}
+\|a\|[\nabla\delta_yv]_{\alpha_0,P_{3\ell}}+[\delta_yj]_{\alpha_0,P_{3\ell}}\nonumber\\
&\lesssim \ell^{-\alpha_0}\|\nabla\delta_yv\|_{\alpha_0,P_{3\ell}}
+[\nabla\delta_yv]_{\alpha_0,P_{3\ell}}+[\delta_yj]_{\alpha_0,P_{3\ell}}.
\end{align}
While obviously
\begin{align}\label{h24}
\|\nabla\delta_yv\|_{P_{3\ell}}\le r^\alpha[\nabla v]_{\alpha,P_{4\ell}},
\end{align}
we need a little argument to see
\begin{align}\label{h25}
[\nabla\delta_yv]_{\alpha_0,P_{3\ell}}+[\delta_yj]_{\alpha_0,P_{3\ell}}
\lesssim r^{\alpha-\alpha_0}\big([\nabla v]_{\alpha,P_{4\ell}}+[j]_{\alpha,P_{4\ell}}\big).
\end{align}
Indeed, let us focus on $j$; given two points $z$, $z'$ in $P_{3\ell}$
we write $\delta_yj(z)-\delta_yj(z')$ in the two ways of $(j(z+(0,y))-j(z))$ $-(j(z'+(0,y))-j(z'))$
and $(j(z+(0,y))-j(z'+(0,y)))$ $-(j(z')-j(z))$ to see that (because of $|y|\le r\le\ell$)
\begin{align*}
|\delta_yj(z)-\delta_yj(z')|\le 2[j]_{\alpha,P_{4\ell}}(\min\{d(z,z'),r\})^\alpha,
\end{align*}
%
%where we call $d(z,z')$ $=\sqrt{|t-t'|}+|x-x'|$,
and thus as desired
\begin{align*}
\frac{|\delta_yj(z)-\delta_yj(z')|}{d^{\alpha_0}(z,z')}
&\le 2[j]_{\alpha,P_{4\ell}}\min\{d^{\alpha-\alpha_0}(z,z'),r^\alpha d^{-\alpha_0}(z,z')\})\nonumber\\
&\le 2[j]_{\alpha,P_{4\ell}} r^{\alpha-\alpha_0}\quad\mbox{since}\;\alpha\ge\alpha_0\ge 0.
\end{align*}
Inserting (\ref{h24}) and (\ref{h25}) into (\ref{h23}) we obtain
\begin{align}\label{h31}
\lefteqn{[a\nabla\delta_yv+\delta_yj]_{\alpha_0,P_{3\ell}}}\nonumber\\
&\lesssim \ell^{-\alpha_0}r^\alpha[\nabla v]_{\alpha,P_{4\ell}}
+r^{\alpha-\alpha_0}\big([\nabla v]_{\alpha,P_{4\ell}}+[j]_{\alpha,P_{4\ell}}\big).
\end{align}
Inserting (\ref{h30}) and (\ref{h31}) into (\ref{h26}) we obtain the iterable form
\begin{align*}
\lefteqn{\inf_c\|\nabla w-c\|_{P_r}}\nonumber\\
&\lesssim\frac{r}{\ell}\inf_c\|\nabla w-c\|_{P_{4\ell}}
+r^{\alpha}(\frac{\ell}{r})^{\alpha_0}\big([\nabla v]_{\alpha,P_{4\ell}}+[j]_{\alpha,P_{4\ell}}\big).
\end{align*}
Relabelling $4\ell$ by $\ell$ we obtain for all $r\le\ell\le L$
\begin{align*}
\lefteqn{r^{-\alpha}\inf_c\|\nabla w-c\|_{P_r}}\nonumber\\
&\lesssim(\frac{r}{\ell})^{1-\alpha}\ell^{-\alpha}\inf_c\|\nabla w-c\|_{P_{\ell}}
+(\frac{\ell}{r})^{\alpha_0}\big([\nabla v]_{\alpha,P_{L}}+[j]_{\alpha,P_{L}}\big).
\end{align*}
By the triangle inequality in $\|\cdot\|$ and by
$\sup_{r\le L}r^{-\alpha}\inf_c\|\nabla v-c\|_{P_r}$ $\le[\nabla v]_{\alpha,P_L}$ this may be upgraded to
\begin{align*}
\lefteqn{r^{-\alpha}\inf_c\|\nabla u-c\|_{P_r}}\nonumber\\
&\lesssim(\frac{r}{\ell})^{1-\alpha}\ell^{-\alpha}\inf_c\|\nabla u-c\|_{P_{\ell}}
+(\frac{\ell}{r})^{\alpha_0}\big([\nabla v]_{\alpha,P_{L}}+[j]_{\alpha,P_{L}}\big).
\end{align*}
Slaving $\ell$ to $r$ via $\ell=Mr$ for some $M\ge1$ to be fixed presently, we obtain
from distinguishing the ranges $r\le\frac{L}{M}$ and $\frac{L}{M}\le r\le L$ that
\begin{align}\label{wi01}
\lefteqn{\sup_{r\le L}r^{-\alpha}\inf_c\|\nabla u-c\|_{P_r}
\lesssim \sup_{\frac{L}{M}\le r\le L}r^{-\alpha}\inf_c\|\nabla u-c\|_{P_r}}\nonumber\\
&+M^{\alpha-1}\sup_{\ell\le L}\ell^{-\alpha}\inf_c\|\nabla u-c\|_{P_\ell}
+M^{\alpha_0}\big([\nabla v]_{\alpha,P_{L}}+[j]_{\alpha,P_{L}}\big).
\end{align}
Clearly, the first right hand side term is controlled as follows
\begin{align*}
\sup_{\frac{L}{M}\le r\le L}r^{-\alpha}\inf_c\|\nabla u-c\|_{P_r}
\le (\frac{M}{L})^{\alpha-\alpha_0}[\nabla u]_{\alpha_0,P_L}
\stackrel{(\ref{h19})}{\le}M^{\alpha-\alpha_0} L^{-\alpha}.
\end{align*}
Hence fixing an $M=M(d,\lambda,\Lambda,\alpha_0,\alpha)$ sufficiently large, we may absorb the
second right hand side term in (\ref{wi01}) into the left hand side to obtain
\begin{align}\label{wi02}
\sup_{r\le L}r^{-\alpha}\inf_c\|\nabla u-c\|_{P_r}
\lesssim L^{-\alpha}+[\nabla v]_{\alpha,P_{L}}+[j]_{\alpha,P_{L}}.
\end{align}
For this, we do not need to know beforehand that the left hand side side is finite, since (\ref{wi01}) also
holds when the two suprema are restricted to $\epsilon\le r\le L$ and $\epsilon\le\ell\le L$ for
any $\epsilon>0$, which is finite since $\nabla u$ is in particular assumed to be continuous.
Hence we obtain (\ref{wi02}) with supremum restricted to $\epsilon\le r\le L$, in which we now
may let $\epsilon\downarrow 0$ to recover the form as stated in (\ref{wi02}).
By the standard norm equivalence
\begin{align*}
\sup_{r\le L}r^{-\alpha}\|\nabla u-\nabla u(0)\|_{P_r}
\lesssim\sup_{r\le L}r^{-\alpha}\inf_c\|\nabla u-c\|_{P_r},
\end{align*}
and shifting the origin into an arbitrary $z\in P_L$, we obtain (\ref{h05}) from (\ref{wi02}).

%%%%%%%%%%%%%%%%%%%%%%%%%%%%%%%%%%%%%%%%%%%%%%%%%%%%%%%%%%%%%%%%%%%%%%%%%%%%

\section{Proof of Corollary~\ref{main}}

In view of Lemmas \ref{L1} and (\ref{L2}) and the definition (\ref{h40}) of the $C^{1,\alpha}$-semi-norm, 
it remains to shows that for all spatial points $x$ and times $t$ and $t'$, the difference $w:=u-v$ satisfies
\begin{align}\label{h41}
|w(t,x)-w(t',x)|\lesssim\big([\nabla u]_\alpha+[j]_\alpha+[\nabla w]_\alpha\big)\sqrt{|t-t'|}^{1+\alpha}.
\end{align}
To this purpose, we rewrite (\ref{h02}) as $\partial_tw$ $=\nabla\cdot(A(\nabla u)+j)$
to which we apply spatial convolution on scale $r$ to be fixed later. This yields the estimate
\begin{align*}
\|\partial_tw_r\|\lessim\frac{1}{r^{1-\alpha}}[A(\nabla u)+j]_{\alpha}
\stackrel{(\ref{h20})}{\lesssim}\frac{1}{r^{1-\alpha}}([\nabla u]_\alpha+[j]_{\alpha}).
\end{align*}
Form this we deduce
\begin{align*}
|w_r(t,x)-w_r(t',x)|\lesssim([\nabla u]_\alpha+[j]_\alpha)\frac{|t-t'|}{r^{1-\alpha}}.
\end{align*}
We may take the convolution kernel $\phi_r$ to be symmetric, so that in particular
$w_r(t,x)=\int\phi_r(x-y)(w(t,y)-\nabla w(t,x)\cdot(y-x))dy$, to the effect of
\begin{align*}
|w(t,x)-w_r(t,x)|\lesssim [\nabla w]_{\alpha}r^{1+\alpha}.
\end{align*}
The last two estimates combine to
\begin{align*}
|w(t,x)-w(t',x)|\lesssim([\nabla u]_\alpha+[j]_\alpha)\frac{|t-t'|}{r^{1-\alpha}}+[\nabla w]_{\alpha}r^{1+\alpha}.
\end{align*}
Optimizing through the choice of $r=\sqrt{t}$ yields (\ref{h41}).

%%%%%%%%%%%%%%%%%%%%%%%%%%%%%%%%%%%%%%%%%%%%%%%%%%%%%%%%%%%%
%%%%%%%    Proof of Lemma  lem:stoch %%%%%%%%%%%%%%%%%%%%%%%%%%%%%%%%%%%%%%

\section{Proof of Lemma~\ref{lem:stoch}}

Throughout the proof we fix $j \in \{1, \ldots, d\}$ and set  $h = \partial_j v$.  We aim to show
that for $C$ large enough  and $\alpha < \min\{ \frac{s - d}{2}, 1 \}$
\begin{align*}
\Big\langle \exp \Big(\frac{1}{C} [h]_\alpha \Big) \Big\rangle < \infty.
\end{align*}
We assume without loss of generality that $\frac{s-d}{2} < 1$. 

\medskip

First we recall that by definition $v$ and $h$ are $1$-periodic in each spatial direction and 
$v(t,x) = h(t,x) =0$ for $t \leq 0$. Furthermore for $t>1$ $h$ solves
\begin{align*}
(\partial_t -\Delta) h =0 
\end{align*}
so that by standard continuity  properties of the heat equation in H\"older norms  we have
$[h]_{\alpha} \lesssim  [h]'_{\alpha}$ where $ [h]'_{\alpha}$ is the local H\"older norm
defined by
\begin{align*}
[h]'_{\alpha}:=\sup_{R\in(0,1)}\frac{1}{R^\alpha}
\sup_{\stackrel{(t,x),(s,y)\in(0,1)\times(-1,1)^d}{\sqrt{|t-s|}+|x-y|<R}}|h(t,x)-h(s,y)|.
\end{align*}
We thus aim to establish
\begin{align}\label{sl1}
\Big\langle \exp \Big(\frac{1}{C} {[h]'_\alpha}^2 \Big) \Big\rangle < \infty.
\end{align}

The core stochastic ingredient for the proof of \eqref{sl1} is the following bound on second moments of 
increments of $h$: For $(t,x)$, $(t',x') \in [0,1] \times \R^d$ we have
\begin{align}\label{aim-first-step}
\big\langle ( h(t,x) - h(t', x') )^2 \big \rangle \lesssim |t-t'|^{\frac{s-d}{2}} + |x-x'|^{s-d}.
\end{align}

\medskip

The argument for \eqref{aim-first-step}  is based on the following Fourier representation for $h$: For $t \in [0,1]$ and $x \in \R^d$
we get by differentiating \eqref{e:def:g} with respect to $x_j$
\begin{align*}
%g(t,x) &= \sum_{k \in (2 \pi \Z)^d}  \sqrt{\hat{K}(k)} e^{i  k \cdot x }\int_0^t e^{-(t-s) |k|^2} d \beta_k(s),\\
h(t,x) &= \sum_{k \in (2 \pi \Z)^d}  \sqrt{\hat{K}(k)}  \; i k_j e^{i k \cdot x }\int_0^t e^{-(t-s) |k|^2} d \beta_k(s),\\
\end{align*}
 which leads to the expression 
 \begin{align}
 \notag
 \big\langle  h(t,x)  h(t', x') \big\rangle &= \sum_{k \in (2 \pi \Z)^d} \hat{K}(k) |k_j|^2 e^{i  k  \cdot (x-x') }\int_0^{t'} e^{-(t-s)|k|^2} e^{- (t'-s)|k|^2} ds   \\
 \label{sp1}
 &= \sum_{k \in (2 \pi \Z)^d} \hat{K}(k)  \frac{|k_j|^2}{2|k|^2} e^{i k  \cdot (x-x') }  \big[   e^{-(t-t') |k|^2} - e^{-(t+t') |k|^2 } \big],
 \end{align}
 valid for $t' \leq t$. In order to deduce \eqref{aim-first-step}, we use the triangle inequality and treat the cases $t =t', \, x \neq x'$ and $t \neq t', \, x = x'$
 separately. In the first case we get
 \begin{align*}
  & \big\langle  ( h(t,x)  -  h(t, x'))^2 \big\rangle \\
  & = 2  \big\langle   h(t,x)^2  -  h(t,x) h(t, x') \big\rangle\\
  & = 2 \sum_{k \in (2 \pi \Z)^d} \hat{K}(k)  \frac{|k_j|^2}{2|k|^2} ( 1-  e^{i k  \cdot (x-x') })  \big[   1 - e^{-(t+t') |k|^2 } \big]. 
 \end{align*}
 Now using the simple estimates $ \frac{|k_j|^2}{2|k|^2} \leq \frac12$, $|  1-  e^{i k  \cdot (x-x') }| \leq \min \{ 2 , | k \cdot (x - x')| \}$
 as well as $ \big[   1 - e^{-(t+t') |k|^2 } \big] \leq 1$, and recalling  condition~\eqref{condition} on $\hat{K}$
  this turns into the estimate
 \begin{align*}
   & \big\langle  ( h(t,x)  -  h(t, x'))^2 \big\rangle \\
  & \lesssim  \sum_{|k| \leq |x - x'|^{-1}}  \hat{K}(k)  |k| |x-x'|   +  \sum_{|k| > |x-x'|^{-1}}  \hat{K}(k) \\
 & \lesssim   |x-x'| \sum_{|k| \leq |x-x'|^{-1}}   \frac{ |k|}{(1 + |k|^2)^\frac{s}{2}}   +  \sum_{|k| > |x-x'|^{-1}}  \frac{ 1}{(1 + |k|^2)^\frac{s}{2}}  \\ 
 & \lessim |x-x'|^{s - d},
 \end{align*}
where $\lesssim$ means $\leq C$.
In the same way we get by specialising \eqref{sp1} to $x = x'$ and treating the case  $t \geq  t'$
 \begin{align*}
  & \big\langle  ( h(t,x)  -  h(t', x))^2 \big\rangle \\
  & =   \big\langle   h(t,x)^2 +  h(t',x)^2   - 2 h(t,x) h(t', x)^2 \big\rangle\\
  & = 2 \sum_{k \in (2 \pi \Z)^d} \hat{K}(k)    \frac{|k_j|^2}{2|k|^2}  \big[   2 - e^{-2t |k|^2 }    - e^{-2t' |k|^2 } - 2    e^{-(t-t') |k|^2} + 2 e^{-(t+t') |k|^2 } \big]. 
 \end{align*}
Now using again $ \frac{|k_j|^2}{2|k|^2}   \leq \frac12$ as well as 
\begin{equation*}
|  2 - e^{-2t |k|^2 }    - e^{-2t' |k|^2 } - 2    e^{-(t-t') |k|^2} + 2 e^{-(t+t') |k|^2 } | \leq 4  \min \{ 1 , |t - t'| |k|^2 \},
\end{equation*}
and using \eqref{condition} once more
this turns into 
\begin{align*}
& \big\langle  ( h(t,x)  -  h(t', x))^2 \big\rangle \\
  & \lessim   |t-t'| \sum_{|k|^2 \leq |t-t'|^{-1}}   \frac{ |k|^2}{(1 + |k|^2)^\frac{s}{2}}   +  \sum_{|k|^2 > |t-t'|^{-1}}  \frac{ 1}{(1 + |k|^2)^\frac{s}{2}}  \\ 
 & \lessim |t-t'|^{\frac{s - d}{2}},
\end{align*}
and thus \eqref{aim-first-step} follows.

\medskip

We now apply Kolmogorov's continuity theorem to $h$; for the convenience of the reader
we give a self-contained argument. We first
appeal to Gaussianity to post-process \eqref{aim-first-step}, which we rewrite as
\begin{equation}\nonumber
\Big\langle\frac{1}{R^{s-d}}(h(t,x)-h(s,y))^2 \Big\rangle \lesssim 1\quad\mbox{provided}\;|t-s|\le 3R^2,|x-y|\le R
\end{equation}
for a given scale $R$. By Gaussianity of $h$ we can upgrade this estimate to 
\begin{align}
\notag
&\Big\langle\exp\Big(\frac{1}{CR^{s-d}}(h(t,x)-h(s,y))^2\Big)\Big\rangle\lesssim 1\\
\label{l3.4}
&\qquad\qquad\qquad\mbox{for}\;|t-s|\le 3R^2,\;|x-y|\le R.
\end{align}

\medskip

Thus  proving the desired estimate \eqref{sl1} on Gaussian moments of the local H\"older-norm 
$[h]'_\alpha$ amounts to exchanging the expectation and the supremum over $(t,x)$, $(s,y)$ in (\ref{l3.4})
at the prize of a decreased H\"older exponent $\alpha<\frac{s-d}{2}$.
To this purpose, we now argue that for $\alpha>0$, the supremum over a continuum can be replaced 
by the supremum over a discrete set: For $R<1$ we define the grid  
\begin{align*}
\Gamma_R= [0,1]\times[-1,1]^d \cap(R^2\mathbb{Z}\times R\mathbb{Z}^d)
\end{align*}
and claim that
\begin{eqnarray}
[h]'_{\alpha}\label{l3.7}
\lesssim\sup_{R}\frac{1}{R^\alpha}
\sup_{\stackrel{(t,x),(s,y)\in\Gamma_R }
{|t-s|\le 3 R^2,|x-y|\le R}}|h(t,x)-h(s,y)|=:\Theta \nonumber,
\end{eqnarray}
where the first $sup$ runs over all $R$ of the form $2^{-N}$ for an integer $N \geq 1$.
Hence we have to show for arbitrary $(t,x),(s,y)\in(-1,0)\times(-1,1)^d$ that
\begin{equation}\label{l3.6}
|h(t,x)-h(s,y)|\lesssim\Theta\big(\sqrt{|t-s|}+|x-y|\big)^\alpha.
\end{equation}
By density, we may assume that $(t,x),(s,y)\in r^2\mathbb{Z}\times r\mathbb{Z}^d$
for some dyadic $r=2^{-N}<1$ (this density argument requires the qualitative a priori
information of the continuity of $h$, which can be circumvented by approximating $h$). 
 For every dyadic level $n=N,N-1,\cdots $ we now recursively construct two sequences $(t_n,x_n)$
$(s_n,y_n)$ of space-time points,
starting from $(t_N,x_N)=(t,x)$ and
$(s_N,y_N)=(s,y)$, with the following properties 
\begin{enumerate}
\item[a)] they are in the corresponding lattice of scale $2^{-n}$, i.\ e.\ we have
$(t_n,x_n),(s_n,x_n)$  $\in (2^{-n})^2\mathbb{Z}\times 2^{-n}\mathbb{Z}^d$,
\item[b)] they are close to their predecessors in the sense of
$|t_{n}-t_{n+1}|,|s_{n}-s_{n+1}|\le 3(2^{-(n+1)})^2$ and  $|x_{n,i}-x_{n+1,i}|,|y_{n,i}-y_{n+1,i}|\le 2^{-(n+1)}$,
where $x_{n,i}$, $x_{n+1,i}$, $\ldots$ denote the $i$-component of $x_{n}$, $x_{n+1}$, $\ldots$.
So  by definition of $\Theta$ we have
\begin{align}
\notag
|h(t_n,x_n)-h(t_{n+1},x_{n+1})| &\lesssim \Theta (2^{-(n+1)})^\alpha,\\
\label{l3.9}
|h(s_n,y_n)-h(s_{n+1},y_{n+1})| 
&\lesssim \Theta (2^{-(n+1)})^\alpha,
\end{align}
and 
\item[c)] such that $|t_n-s_n|$ and $|x_n-y_n|$ are minimized among the points satisfying a) and b).
\end{enumerate}
Because of the latter, we have
\begin{equation}\nonumber
(t_M,x_M)=(s_M,y_M)\quad\mbox{for some}\;M\;\mbox{with}\quad 2^{-M}\le\max\{\sqrt{|t-s|},|x-y|\},
\end{equation}
so that by the triangle inequality we gather from (\ref{l3.9})
\begin{equation}\nonumber
|h(t,x)-h(s,y)|\lesssim \sum_{n=N-1}^M\Theta (2^{-(n+1)})^\alpha\le\Theta\frac{(2^{-M})^\alpha}{2^\alpha-1},
\end{equation}
which yields (\ref{l3.6}).

\medskip
Equipped with (\ref{l3.7}), we now may upgrade (\ref{l3.4}) to \eqref{sl1}.
Indeed, (\ref{l3.7}) can be reformulated
on the level of characteristic functions as
\begin{equation}\nonumber
I\big(([h]'_{\alpha})^2\ge M)
\le\sup_{R }\max_{(t,x),(s,y) \in \Gamma_R}I\Big(\frac{1}{R^{s-d}}(h(t,x)-h(s,y))^2\ge\frac{M}{CR^{s-d-2\alpha}}\Big),
\end{equation}
where as in (\ref{l3.7}) $R$ runs over all $2^{-N}$ for integers $N \geq 1$. Replacing the suprema by sums
in order to take the expectation, we obtain
\begin{align*}
&\big\langle I\big(([h]'_{\alpha})^2\ge M)\big\rangle \\
&\le\sum_{R}\sum_{(t,x),(s,y)} \Big\langle I\Big(\frac{1}{R^{s-d}}(h(t,x)-h(s,y))^2\ge\frac{M}{CR^{s-d-2\alpha}}\Big)\Big\rangle. 
\end{align*}
We now appeal to Chebyshev's inequality in order to make use of (\ref{l3.4}):
\begin{eqnarray*}\nonumber
&& \big\langle I\big( ([h]'_{\alpha})^2\ge M)\big\rangle\\
%&\lesssim&\sum_{R}\sum_{(t,x),(s,y)}   \\
&\lesssim&\sum_{R}\sum_{(t,x),(s,y)}\exp \Big(-\frac{M}{CR^{s-d-2\alpha}} \Big)\\
&\lesssim&\sum_{R}\frac{1}{R^{2+d}}\exp\Big(-\frac{M}{CR^{s-d-2\alpha}}\Big)\\
&\stackrel{R\le 1,M\ge 1}{\le}
&\exp(-\frac{M}{C})\sum_{R}\frac{1}{R^{2+d}}\exp(-\frac{1}{C}(\frac{1}{R^{s-d-2\alpha}}-1))
\lesssim \exp(-\frac{M}{C}),
\end{eqnarray*}
where in the second step we have used that the number of pairs $(t,x),(s,y)$
of neighboring lattice points is bounded by $C\frac{1}{R^{2+d}}$ and in the last
step we have used that stretched exponential decay (recall $s-d-2\alpha>0$) beats polynomial growth.
The last estimate immediately yields (\ref{sl1}).

\bibliographystyle{abbrv}
\bibliography{HoelderDivergence}

\end{document}